\documentclass[a4paper,12pt]{article}

\usepackage{float}
\usepackage{amsmath} 
\usepackage{amssymb}
\usepackage{amsfonts}
\usepackage{epstopdf}
\usepackage{color,fancybox,graphicx}
\usepackage{url}
\usepackage{psfrag}

\numberwithin{equation}{section}

\newtheorem{lem}{Lemma}[section]
\newtheorem{cor}{Corollary}[section]
\newtheorem{pro}{Proposition}[section]
\newtheorem{theo}{Theorem}[section]

\newcommand{\bX}{\mathbf{X}}
\newcommand{\bx}{\mathbf{x}}

\begin{document}

\begin{center}

{\sc \Large An Affine Invariant $k$-Nearest Neighbor 
\medskip

Regression Estimate}
\medskip

\end{center}
{\bf G\'erard Biau\footnote{Corresponding author.}\\
{\it Universit\'e Pierre et Marie Curie\footnote{Research partially supported by the French National Research Agency under grant ANR-09-BLAN-0051-02 ``CLARA''.} \& Ecole Normale Sup{\'e}rieure\footnote{Research carried out within the INRIA project ``CLASSIC'' hosted by Ecole Normale Sup{\'e}rieure and CNRS.}, France}}\\
\textsf{gerard.biau@upmc.fr}
\bigskip

{\bf Luc Devroye}\\
{\it McGill University, Canada}\footnote{Research sponsored by NSERC Grant A3456 and FQRNT Grant 90-ER-0291.}\\
\textsf{lucdevroye@gmail.com}
\bigskip

{\bf Vida Dujmovi\'c}\\
{\it Carleton University, Canada}\footnote{Research sponsored by NSERC Grant RGPIN 402438-2011.}\\
\textsf{vida@scs.carleton.ca}
\bigskip

{\bf Adam Krzy\.zak}\\
{\it Concordia University, Canada}\footnote{Research sponsored by NSERC Grant N00118.}\\
\textsf{krzyzak@cs.concordia.ca}
\medskip

\begin{abstract}
\noindent {\rm We design a data-dependent metric in $\mathbb R^d$ and use it
to define the $k$-nearest neighbors of a given point.
Our metric is invariant under all affine transformations.
We show that, with this metric, the standard $k$-nearest neighbor regression
estimate is asymptotically consistent under the usual
conditions on $k$, and minimal requirements on the
input data.
\medskip

\noindent \emph{Index Terms} --- Nonparametric estimation,
Regression function estimation,
Affine invariance,
Nearest neighbor methods,
Mathematical statistics.
\medskip

\noindent \emph{2010 Mathematics Subject Classification}: 62G08, 62G05, 62G20.}

\end{abstract}

\section{Introduction}
The prediction error of standard nonparametric regression methods may be critically affected by a linear transformation of the coordinate axes.  It is typically the case for the popular $k$-nearest neighbor ($k$-NN) predictor (Fix and Hodges \cite{Fix1, Fix2}, Cover and Hart \cite{coverandhart}, Cover \cite{Cover68a, Cover68b}), where a mere rescaling of the coordinate axes has a serious impact on the capabilities of this estimate. This is clearly an undesirable feature, especially in applications where the data measurements represent physically different quantities, such as temperature, blood pressure, cholesterol level, and the age of the patient. In this example, a simple change in, say, the unit measure of the temperature parameter will lead to totally different results, and will thus force the statistician to use a somewhat arbitrary preprocessing step prior to the $k$-NN estimation process.  Furthermore, in several practical implementations, one would like, for physical or economical reasons, to supply the freshly collected data to some machine without preprocessing. 
\medskip

In this paper, we discuss a variation of the $k$-NN regression estimate whose definition is not affected by affine transformations of the coordinate axes. Such a modification could save the user a subjective preprocessing step and would save the manufacturer the trouble of adding input specifications.
\medskip

The data set we have collected can be regarded as a collection of independent and identically distributed $\mathbb R^d\times \mathbb R$-valued random variables $\mathcal D_n=\{(\bX_1,Y_1),\hdots,(\bX_n,Y_n)\}$, independent of and with the same distribution as a generic pair $(\bX, Y)$ satisfying $\mathbb E |Y| < \infty$. The space $\mathbb R^d$ is equipped with the standard Euclidean norm $\|.\|$. For  fixed $\bx \in \mathbb R^d$, our goal is to estimate the regression function $r(\bx)=\mathbb E[Y|\bX=\bx]$ using the data $\mathcal D_n$. In this context, the usual $k$-NN regression estimate takes the form
$$r_n(\bx;\mathcal D_n)=\frac{1}{k_n}\sum_{i=1}^{k_n}Y_{(i)}(\bx),$$
where $(\bX_{(1)}(\bx),Y_{(1)}(\bx)),\hdots,(\bX_{(n)}(\bx),Y_{(n)}(\bx))$ is a reordering of the data according to increasing distances $\|\bX_i-\bx\|$ of the $\bX_i$'s to $\bx$. (If distance ties occur, a tie-breaking strategy must be defined. For example, if $\|\bX_i-\bx\|=\|\bX_j-\bx\|$, $\bX_i$ may be declared ``closer'' if $i <j$, i.e., the tie-breaking is done by indices.)  For simplicity, we will suppress $\mathcal D_n$ in the notation and write $r_n(\bx)$ instead of $r_n(\bx;\mathcal D_n)$. Stone \cite{Stone} showed that, for all $p\geq 1$, $\mathbb E [r_n(\bX)-r(\bX)]^p \to 0$ for all possible distributions of $(\bX,Y)$ with $\mathbb E |Y|^p <\infty$, whenever $k_n \to \infty$ and $k_n/n \to 0$ as $n \to \infty$. Thus, the $k$-NN estimate behaves asymptotically well, without exceptions. This property is called $L_p$ universal consistency.
\medskip

Clearly, any affine transformation of the coordinate axes influences the $k$-NN estimate through the norm $\|.\|$, thereby illuminating an unpleasant face of the procedure. To illustrate this remark, assume that a nontrivial affine transformation $T: \mathbf z \mapsto A\mathbf z +\mathbf b$ (that is, a nonsingular linear transformation $A$ followed by a translation $\mathbf b$) is applied to both $\bx$ and $\bX_1,\hdots,\bX_n$. Examples include any number of combinations of rotations, translations, and linear rescalings. Denote by $\mathcal D'_n=(T(\bX_1),Y_1),\hdots,(T(\bX_n),Y_n)$ the transformed sample. Then, for such a function $T$, one has $r_n(\bx;\mathcal D_n)\neq r_n(T(\bx);\mathcal D'_n)$ in general, whereas $r(\bX)=\mathbb E[Y|T(\bX)]$ since $T$ is bijective. Thus, to continue our discussion, we are looking in essence for a regression estimate $r_n$ with the following property:
\begin{equation}
\label{neige}
r_n(\bx;\mathcal D_n)=r_n(T(\bx);\mathcal D'_n).
\end{equation}
We call $r_n$ affine invariant. Affine invariance is indeed a very strong but highly desirable property. In $\mathbb R^d$, in the context of $k$-NN estimates, it suffices to be able to define an affine invariant distance measure, which is necessarily data-dependent. With this objective in mind, we develop in the next section an estimation procedure featuring (\ref{neige}) which in form coincides with the $k$-NN estimate, and establish its consistency in Section 3. Proofs of the most technical results are gathered in Section 4.
\medskip

It should be stressed that what we are after in this article is an estimate of $r$ which is invariant by an affine transformation
of {\it both} the query point $\bx$ and the original regressors $\bX_1, \hdots, \bX_n$. When the sole regressors are subject to 
such a transformation, it is then more natural to talk of ``affine equivariant'' regression estimates rather than of ``affine invariant" ones; this is more in line with the terminology used, for example, in Ollila, Hettmansperger, and Oja \cite{Ollila1} and Ollila, Oja, and Koivunen \cite{Ollila2}. These affine invariance and affine equivariance requirements, however, are strictly equivalent.
\medskip

There have been many attempts in the nonparametric literature
to achieve affine invariance. One of the most natural ones relates to the so-called transformation-retransformation proposed by Chakraborty, Chaud\-huri, and Oja \cite{Chakra}.  
That method and many variants have been discussed 
in texts such as \cite{DGL} and \cite{Laciregressionbook} 
for pattern recognition and regression, respectively,
but they have also been used in kernel density estimation (see, e.g.,
Samanta \cite{Samanta}). It is worth noting that, computational issues aside, the
transformation step (i.e., premultiplication of the regressors by $\hat{M}_n^{-1}$, where $\hat{M}_n$ is an affine equivariant
scatter estimate) may be based on a statistic $\hat{M}_n$ that does not require finiteness of any
moment. A typical example is the scatter estimate proposed in Tyler \cite{Tyler} or Hettmansperger and Randles \cite{Hettman2}. 
Rather, our procedure
takes ideas from the classical nonparametric literature using
concepts such as multivariate ranks. It is closed in spirit of the approach of Paindaveine and Van Bever \cite{Davy}, who introduce a class of depth-based classification procedures that are of a nearest neighbor nature.
\medskip

There are also attempts at getting invariance to other
transformations. The most important concept here is that of invariance
under monotone transformations of the coordinate axes.  In
particular, any strategy that uses only the coordinatewise ranks
of the $\bX_i$'s achieves this.  The onus, then, is to show
consistency of the methods under the most general conditions
possible.  For example, using an $L_p$ norm on the $d$-vectors of differences
between ranks, one can show that the classical $k$-NN regression function estimate is universally consistent in the sense of Stone \cite{Stone}.
This was observed by Olshen \cite{Olshen}, and shown by Devroye \cite{Devroye4}
(see also Gordon and Olshen \cite{Gordon2,Gordon1}, Devroye and Krzy\.zak \cite{DK}, and Biau and Devroye \cite{BD} for related works).
Rules based upon statistically equivalent blocks (see, e.g., Anderson \cite{Anderson}, Quesenberry and Gessaman \cite{Quesenberry}, Gessaman \cite{Gessaman1}, Gessaman and Gessaman \cite{Gessaman2}, and Devroye, Gy\"orfi, and Lugosi \cite[Section 21.4]{DGL}) are other important examples of regression methods invariant with respect to monotone transformations of the coordinate axes. These methods and their generalizations partition the space with sets that contain a fixed number of data points each.
\medskip

It would be interesting to consider in a future paper the possibility
of morphing the input space in more general ways than those suggested in
the previous few paragraphs of the present article.  It should be possible,
in principle, to define appropriate metrics to obtain invariance
for interesting large classes of nonlinear transformations,
and show consistent asymptotic behaviors.
\section{An affine invariant $k$-NN estimate}
The $k$-NN estimate we are discussing is based upon the notion of empirical distance. Throughout, we assume that the distribution of $\bX$ is absolutely continuous with respect to the Lebesgue measure on $\mathbb R^d$ and that $n\geq d$. Because of this density assumption, any collection $\bX_{i_1},\hdots,\bX_{i_d}$ ($1\leq i_1<i_2<\hdots<i_d\leq n)$ of $d$ points among $\bX_1, \hdots, \bX_n$ are in general position with probability 1. Consequently, there exists with probability 1 a unique hyperplane in $\mathbb R^d$ containing these $d$  random points, and we denote it by $\mathcal H(\bX_{i_1},\hdots,\bX_{i_d})$.  
\medskip

\begin{figure}[h]
\psfrag{0}{$0$}
\psfrag{1}{$1$}
\psfrag{2}{$2$}
\psfrag{3}{$3$}
\psfrag{4}{$4$}
\psfrag{5}{$5$}
\psfrag{p}{$\bx'$}
\psfrag{x}{$\bx$}
\centering
\includegraphics*[height=7cm]{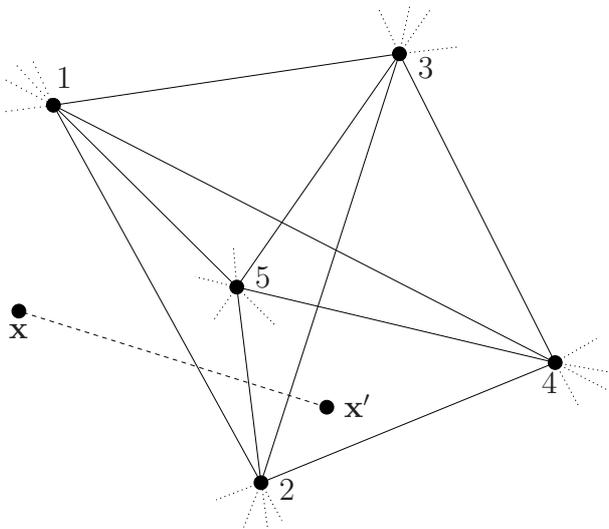}
\caption{An example in dimension 2. The empirical distance between $\bx$ and $\bx'$ is 4. (Note that the hyperplane defined by the pair $(3,5)$ indeed cuts the segment $(\bx,\bx')$, so that the distance is 4, not 3.)}
\label{figure0}
\end{figure}

With this notation, the empirical distance between $d$-vectors $\bx$ and $\bx'$ is defined as
$${\small \rho_n(\bx,\bx')=\sum_{{1\leq i_1<\hdots <i_d\leq n}}\mathbf 1_{\{\mbox{segment } (\bx,\bx')\mbox{ intersects the hyperplane }\mathcal H(\bX_{i_1},\hdots,\bX_{i_d})\}}.}$$
Put differently, $\rho_n(\bx,\bx')$ just counts the number of hyperplanes in $\mathbb R^d$ passing through $d$ out of the points $\bX_1,\hdots,\bX_n$,  that are separating $\bx$ and $\bx'$. Roughly, ``near'' points have fewer intersections, see Figure \ref{figure0} that depicts an example in dimension $2$.
\medskip

This hyperplane-based concept of distance is known in the multivariate rank tests literature as the empirical lift-interdirection function (Oja and Paindaveine \cite{OP}, see also Randles \cite{Randles}, Oja \cite{Oja}, and Hallin and Paindaveine \cite{HP} for companion concepts). It was originally mentioned (but not analyzed) in Hettmansperger, M\"ott\"onen, and Oja \cite{Hett}, and independently suggested as an affine invariant alternative to ordinary metrics in the monograph of Devroye, G\"yorfi, and Lugosi \cite[Section 11.6]{DGL}. We speak throughout of distance even though, for a fixed sample of size $n$, $\rho_n$ is only defined with probability 1 and is not a distance measure {\it stricto sensu} (in particular, $\rho_n(\bx, \bx') = 0$ does
not imply that $\bx = \bx'$). Nevertheless, this empirical distance is invariant under affine transformations $\bx \mapsto A\bx+\mathbf b$, where $A$ is some arbitrary nonsingular linear map and $\mathbf b$ any offset vector (see, for instance, Oja and Paindaveine \cite[Section 2.4]{OP}).
\medskip

Now, fix $\bx \in \mathbb R^d$ and let $\rho_n(\bx, \bX_i)$ be the empirical distance between $\bx$ and some observation $\bX_i$ in the sample $\bX_1, \hdots, \bX_n$. (That is, the number of hyperplanes in $\mathbb R^d$ passing through $d$ out of the observations $\bX_1, \hdots, \bX_n$, that are cutting the segment $(\bx,\bX_i)$). In this context, the $k$-NN estimate we are considering still takes the familiar form 
\begin{equation*}
r_n(\bx)=\frac{1}{k_n}\sum_{i=1}^{k_n}Y_{(i)}(\bx),
\end{equation*}
with the important difference that now the data set $(\bX_1,Y_1),\hdots,(\bX_n,Y_n)$ is reordered according to increasing values of the empirical distances $\rho_n(\bx,\bX_i)$, {\it not} the original Euclidean metric. By construction, the estimate $r_n$ has the desired affine invariance property and, moreover, it coincides with the standard (Euclidean) estimate in dimension $d=1$.
In the next section, we prove the following theorem.  The distribution of the random variable $\bX$ is denoted by $\mu$.
\begin{theo}[Pointwise $L_p$ consistency]
\label{theoreme1}
\label{convergenceponctuelle}
Assume that $\bX$ has a probability density, that $Y$ is bounded, and that the regression function $r$ is  $\mu$-almost surely continuous.  Then, for $\mu$-almost all $\bx\in \mathbb R^d$ and all $p\geq 1$, if $k_n \to \infty$ and $k_n/n \to 0$,
$$\mathbb E \left |r_n(\bx)-r(
\bx)\right|^p \to 0\quad \mbox{as } n\to \infty.$$
\end{theo}
The following corollary is a consequence of Theorem \ref{convergenceponctuelle} and the Lebesgue dominated convergence theorem.
\begin{cor}[Global $L_p$ consistency]
\label{convergenceglobale}
Assume that $\bX$ has a probability den\-si\-ty, that $Y$ is bounded, and that the regression function $r$ is $\mu$-almost surely continuous.  Then, for all $p\geq 1$, if $k_n \to \infty$ and $k_n/n \to 0$,
$$ \mathbb E\left |r_n(\bX)-r(
\bX)\right|^p \to 0\quad \mbox{as } n\to \infty.$$
\end{cor}
The conditions of Stone's universal consistency theorem given in \cite{Stone}
are not fulfilled for our estimate. For the standard nearest neighbor
estimate, a key result used in the consistency proof by Stone is that
a given data point cannot be the nearest neighbor of more than a constant
number (say, $3^d$) other points.  Such a universal constant
does not exist after our transformation is applied. That means that 
a single data point can have a large influence on the regression
function estimate.  While this by itself does not imply that the 
estimate is not universally consistent, it certainly indicates that any such
proof will require new insights. The addition of two smoothness
constraints, namely that $\bX$ has a density (without, however, imposing
any continuity conditions on the density itself) and that $r$ is $\mu$-almost
surely continuous, is sufficient.
\medskip

The complexity of our procedure in terms of sample size $n$ and
dimension $d$ is quite high. There are ${n \choose d}$ possible
choices of hyperplanes through $d$ points. 
This collection of hyperplanes defines an arrangement, or partition of $\mathbb R^d$
into polytopal regions, also called cells or chambers.
Within each region, the distance to each data point is constant,
and thus, a preprocessing step might consist of setting up a
data structure for determining to which cell a given point $\bx \in \mathbb R^d$
belongs: This is called the point location problem.
Meiser \cite{Meiser} showed that such a data structure exists with
the following properties: $(1)$ it takes space $\mathcal O(n^{d+\varepsilon})$ for any fixed
$\varepsilon > 0$, and $(2)$ point location can be performed in $\mathcal O(\log n)$ time.
Chazelle's cuttings \cite{Chazelle} improve $(1)$ to $\mathcal O(n^d)$. 
Chazelle's processing time for setting up the data structure is $\mathcal O(n^d)$.
Still in the preprocessing step, one can determine for each cell
in the arrangement the distances to all $n$ data points: This can
be done by walking across the graph of cells or by brute force. When
done naively, the overall set-up complexity is $\mathcal O(n^{2d+1})$.
For each cell, one might keep a pointer to the $k$ nearest neighbors.
Therefore, once set up, the computation of the regression function
estimate takes merely $\mathcal O(\log n)$ time for point location, and $\mathcal O(k)$ time
for retrieving the $k$ nearest neighbors.
\medskip

One could envisage a reduction in the complexity by defining
the distances not in terms of all hyperplanes that cut a line segment,
but in terms of the number of randomly drawn hyperplanes that make
such a cut, where the number of random draws is now a carefully
selected number.  By the concentration of binomial random
variables, such random estimates of the distances are expected to
work well, while keeping the complexity reasonable. This idea will
be explored elsewhere.
\section{Proof of the theorem}
Recall, since $\bX$ has a probability density with respect to the Lebesgue measure on $\mathbb R^d$, that any collection $\bX_{i_1},\hdots,\bX_{i_d}$ ($1\leq i_1<i_2<\hdots<i_d\leq n)$ of $d$ points among $\bX_1, \hdots, \bX_n$ defines with probability 1 a unique hyperplane $\mathcal H(\bX_{i_1},\hdots,\bX_{i_d})$ in $\mathbb R^d$. Thus, in the sequel, since no confusion is possible, we will freely refer to ``the hyperplane $\mathcal H(\bX_{i_1},\hdots,\bX_{i_d})$ defined by $\bX_{i_1},\hdots,\bX_{i_d}$'' without further explicit mention of the probability 1 event.
\medskip

Let us first fix some useful notation. The distribution of the random variable $\bX$ is denoted by $\mu$ and its density with respect to the Lebesgue measure is denoted by $f$. For every $\varepsilon >0$, we let $\mathcal B_{\bx,\varepsilon}=\{\mathbf y \in \mathbb R^d:\|\mathbf y-\bx\|\leq \varepsilon\}$ be the closed Euclidean ball with center at $\bx$ and radius $\varepsilon$. We write $A^c$ for the complement of a subset $A$ of $\mathbb R^d$. For two random variables $Z_1$ and $Z_2$, the notation
$$Z_1\leq_{\mbox{\footnotesize st}}Z_2$$
means that $Z_1$ is stochastically dominated by $Z_2$, that is, for all $t\in \mathbb R$,
$$\mathbb P  \{Z_1> t\} \leq \mathbb P\{Z_2> t\}.$$ 

\begin{figure}[h]
\psfrag{x}{$\bx$}
\psfrag{theta}{${\Theta}$}
\psfrag{epsilon}{$\frac{\varepsilon}{2}$}
\psfrag{N1}{$N^2_{\bx, \varepsilon}({\Theta})$}
\psfrag{N2}{$N^1_{\bx, \varepsilon}({\Theta})$}
\psfrag{Q1}{$\mathcal R^2_{\bx, \varepsilon}({\Theta})$}
\psfrag{Q2}{$\mathcal R^1_{\bx, \varepsilon}({\Theta})$}
\psfrag{S}{$\mathcal S_{\bx, \varepsilon}({\Theta})$}
\psfrag{L}{$\mathcal L_{\bx}({\Theta})$}
\centering
\includegraphics*[height=7cm]{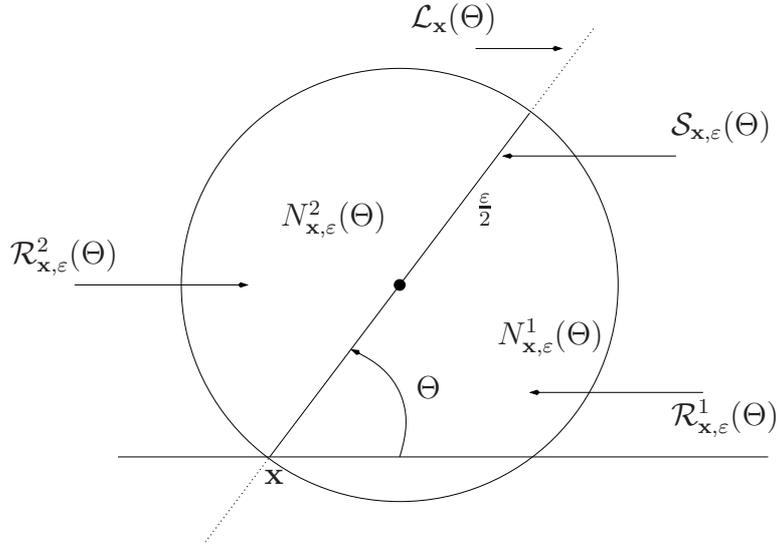}
\caption{The ball $\mathcal B_{\bx,\varepsilon}(\Theta)$ and related notation. Illustration in dimension $2$.}
\label{angle}
\end{figure}

Our first goal is to show that for $\mu$-almost all $\bx$, as $k_n/n\to0$, the quantity $\max_{i=1,\hdots,k_n}\| \bX_{(i)}(\bx)-\bx\|$ converges to $0$ in probability, i.e., for every $\varepsilon >0$,
\begin{equation}
\lim_{n\to\infty}\mathbb P \left\{\max_{i=1,\hdots,k_n}\| \bX_{(i)}(\bx)-\bx\|>\varepsilon\right\}=0.
\label{convergenceenproba}
\end{equation}
So, fix such a positive $\varepsilon$. Let $\delta$ be a real number in $(0,\varepsilon)$ and $\gamma_n$ be a positive real number (eventually function of $\bx$ and $\varepsilon$) to be determined later. To prove identity (\ref{convergenceenproba}), we use the following decomposition, which is valid for all $\bx\in \mathbb R^d$:
\begin{align}
& \nonumber \mathbb P \left\{ \max_{i=1,\hdots,k_n}\| \bX_{(i)}(\bx)-\bx\| > \varepsilon\right\} \\
& \nonumber \quad \leq \mathbb P \left\{ \min_{i=1,\hdots,n \atop \bX_i \in \mathcal B^c_{\bx,\varepsilon}}\rho_n(\bx,\bX_i) <\gamma_n\right\}\\
&\nonumber \qquad +\mathbb P \left\{ \max_{i=1,\hdots,n \atop \bX_i \in \mathcal B_{\bx,\delta}}\rho_n(\bx,\bX_i) \geq \gamma_n\right\}\\
&\nonumber \qquad +\mathbb P \left\{\mbox{Card}\left\{i=1,\hdots,n:\|\bX_i-\bx\| \leq \delta\right\}<k_n\right\}\\
& \quad := {\bf A} + {\bf B} + {\bf C}.\label{grosseinegalite}
\end{align}
The convergence to $0$ of each of the three terms above---from which identity (\ref{convergenceenproba}) immediately follows---are separately analyzed in the next three paragraphs.
\paragraph{Analysis of ${\bf A}$.}  As for now, taking an affine geometry point of view, we keep $\bx$ fixed and see it as the origin of the space. Recall that each point in the Euclidean space $\mathbb R^d$ (with the origin at $\bx$) may be described by its hyperspherical coordinates (see, e.g., Miller \cite[Chapter 1]{Miller}), which consist of a nonnegative radial coordinate $r$ and $d-1$ angular coordinates $\theta_1, \hdots, \theta_{d-1}$, where $\theta_{d-1}$ ranges over $[0,2\pi)$ and the other angles range over $[0,\pi]$ (adaptation of this definition to the cases $d=1$ and $d=2$ is clear). For a $(d-1)$-dimensional vector ${\Theta}=(\theta_1, \hdots, \theta_{d-1})$ of hyperspherical angles, we let $\mathcal B_{\bx, \varepsilon}({\Theta})$ be the unique closed ball anchored at $\bx$ in the direction ${\Theta}$ and with diameter $\varepsilon$ (see Figure \ref{angle} which depicts an illustration in dimension $2$). We also let $\mathcal L_{\bx}({\Theta})$ be the axe defined by $\bx$ and the direction $\Theta$, and let as well $\mathcal S_{\bx, \varepsilon}({\Theta})$ be the open segment obtained as the intersection of $\mathcal L_{\bx}({\Theta})$ and the interior of $\mathcal B_{\bx, \varepsilon}({\Theta})$.
\medskip

Next, for fixed $\bx$, $\varepsilon$ and ${\Theta}$, we split the ball $\mathcal B_{\bx, \varepsilon}({\Theta})$ into $2^{d-1}$ disjoint regions $\mathcal R^1_{\bx, \varepsilon}({\Theta}), \hdots, \mathcal R^{2^{d-1}}_{\bx, \varepsilon}({\Theta})$ as follows. First, the Euclidean space $\mathbb R^d$ is sequentially divided into $2^{d-1}$ symmetric quadrants rotating around the axe $\mathcal L_{\bx}(\Theta)$ (boundary equalities are broken arbitrarily). Next, each region $\mathcal R^j_{\bx, \varepsilon}(\Theta)$ is obtained as the intersection of one of the $2^{d-1}$ quadrants and the ball  $\mathcal B_{\bx, \varepsilon}({\Theta})$. 
\medskip

The numbers of sample points falling in each of these regions are denoted hereafter by $N^1_{\bx,\varepsilon}({\Theta}), \hdots, N^{2^{d-1}}_{\bx,\varepsilon}({\Theta})$ (see Figure \ref{angle}). Letting finally $V_d$ be the volume of the unit $d$-dimensional Euclidean ball, we are now in a position to control the first term of inequality (\ref{grosseinegalite}).
\begin{pro}
\label{ulm1}
For $\mu$-almost all $\bx \in \mathbb R^d$ and all $\varepsilon>0$ small enough, 
$$\mathbb P \left\{ \min_{i=1,\hdots,n \atop \bX_i \in \mathcal B^c_{\bx,\varepsilon}}\rho_n(\bx,\bX_i) <\gamma_n\right\} \to 0 \quad \mbox{as } n \to \infty,$$
provided
$$\gamma_n=n^d\left( \frac{V_d }{2^{2d+1}}\varepsilon^d f(\bx)\right)^{2^{d-1}}.$$
\end{pro}
{\bf Proof of Proposition \ref{ulm1}}\quad Set
$$p_{\bx, \varepsilon}=\min_{j=1, \hdots, 2^{d-1}}\inf_{{\Theta}} \mu \left\{ \mathcal R^j_{\bx,\varepsilon}({\Theta})\right\},$$
where the infimum is taken over all possible hyperspherical angles ${\Theta}$. We know, according to technical Lemma \ref{samedi26}, that for $\mu$-almost all $\bx$ and all $\varepsilon>0$ small enough,
 \begin{equation}
 \label{chypre}
 p_{\bx, \varepsilon}\geq \frac{V_d}{2^{2d}}{\varepsilon}^df(\bx)>0.
 \end{equation}
Thus, in the rest of the proof, we fix such an $\bx$ and assume that $\varepsilon$ is small enough so that the inequalities above are satisfied.
\medskip

\begin{figure}[h]
\psfrag{x}{$\bx$}
\psfrag{X1}{$\bX$}
\psfrag{X2}{$\bX^{\star}$}
\psfrag{e}{$\varepsilon$}
\psfrag{N}{$N^2_{\bx,\varepsilon}(\Theta^{\star})$}
\psfrag{S}{$N^1_{\bx,\varepsilon}(\Theta^{\star})$}
\psfrag{t}{${\Theta}^{\star}$}
\centering
\includegraphics*[height=8cm]{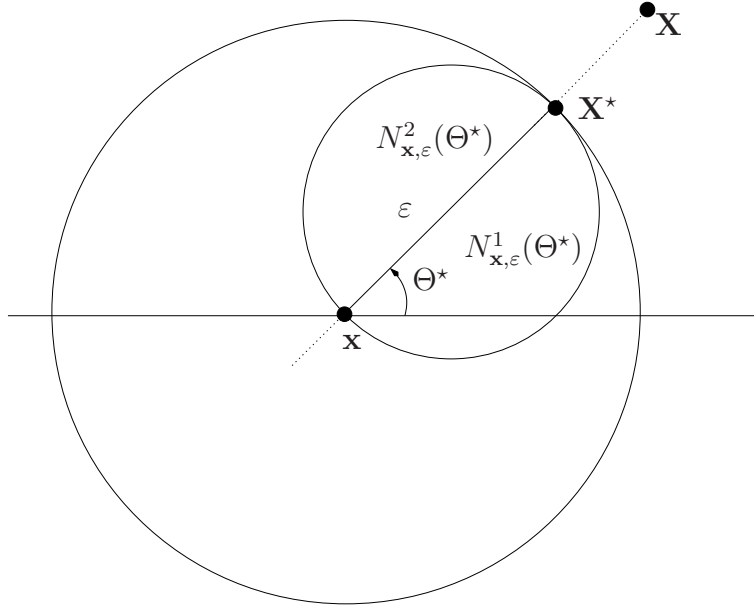}
\caption{The ball $\mathcal B_{\bx, \varepsilon}({\Theta}^{\star})$ in dimension 2.}
\label{anglebis}
\end{figure}

Let $\bX^{\star}$ be defined as the intersection of the line $(\bx,\bX)$ with $\mathcal B_{\bx,\varepsilon}$, and let ${\Theta}^{\star}$ be the (random) hyperspherical angle corresponding to $\bX^{\star}$ (see Figure \ref{anglebis} for an example in dimension 2).
\medskip

Denote by $N_{\bx,\varepsilon}({\Theta^\star})$ the number of hyperplanes passing through $d$ out of the observations $\bX_1,\hdots,\bX_n$ and cutting the segment $\mathcal S_{\bx,\varepsilon}({\Theta^\star})$.  We have
\begin{align*}
\mathbb P \left\{ \min_{i=1,\hdots,n \atop \bX_i \in \mathcal B^c_{\bx,\varepsilon}}\rho_n(\bx,\bX_i) <\gamma_n\right\}&\leq n \mathbb P \left \{\rho_n(\bx,\bX^{\star})<\gamma_n\right\}\\
& =n \mathbb P \left \{N_{\bx,\varepsilon}(\Theta^{\star})<\gamma_n\right\}\\
& \leq n \mathbb P \left \{ \frac{N^1_{\bx,\varepsilon}({\Theta}^{\star})\hdots N^{2^{d-1}}_{\bx, \varepsilon}({\Theta}^{\star})}{n^{2^{d-1}-d}}<\gamma_n\right \},
\end{align*}
where the last inequality follows from technical Lemma \ref{lundi28}. Thus,
\begin{align*}
\mathbb P \left\{ \min_{i=1,\hdots,n \atop \bX_i \in \mathcal B^c_{\bx,\varepsilon}}\rho_n(\bx,\bX_i) <\gamma_n\right\} &\leq n\sum_{j=1}^{2^{d-1}} \mathbb P \left \{ N^j_{\bx, \varepsilon}({\Theta}^{\star}) < \left(\gamma_nn^{2^{d-1}-d}\right)^{1/2^{d-1}}\right\}\\
& = n\sum_{j=1}^{2^{d-1}} \mathbb P \left \{ N^j_{\bx, \varepsilon}({\Theta}^{\star}) < \gamma_n^{1/2^{d-1}}n^{1-d/2^{d-1}}\right\}.
\end{align*}
Clearly, conditionally on $\Theta^{\star}$, each $N^j_{\bx,\varepsilon}({\Theta}^{\star})$ satisfies
$$\mbox{Binomial }(n,p_{\bx,\varepsilon}) \leq_{\mbox{\footnotesize st}} N^j_{\bx,\varepsilon}({\Theta}^{\star})$$
and consequently, by inequality (\ref{chypre}),
$$\mbox{Binomial }\left(n,\frac{V_d}{2^{2d}}\varepsilon^d f(\bx)\right) \leq_{\mbox{\footnotesize st}} N^j_{\bx,\varepsilon}({\Theta}^{\star}).$$
Thus, for each $j=1, \hdots, 2^{d-1}$, by Hoeffding's inequality for binomial random variables (Hoeffding \cite{Hoeffding}), we are led to
\begin{align*}
&\mathbb P \left \{ N^j_{\bx, \varepsilon}({\Theta}^{\star}) < \gamma_n^{1/2^{d-1}}n^{1-d/2^{d-1}}\right\} \\
& \quad= \mathbb E \left [ \mathbb P \left \{ N^j_{\bx, \varepsilon}({\Theta}^{\star}) < \gamma_n^{1/2^{d-1}}n^{1-d/2^{d-1}} | \Theta^{\star}\right\}\right] \\
& \quad \leq \exp \left[-2\left( \gamma_n^{1/2^{d-1}}n^{1-d/2^{d-1}}-n\frac{V_d}{2^{2d}}\varepsilon^d f(\bx)\right)^2/n\right]
\end{align*}
as soon as $\gamma_n^{1/2^{d-1}}n^{1-d/2^{d-1}} < n\frac{V_d}{2^{2d}}\varepsilon^d f(\bx)$. Therefore, taking
$$\gamma_n=n^d\left(\frac{V_d}{2^{2d+1}}\varepsilon^d f(\bx)\right)^{2^{d-1}},$$
we obtain 
$$\mathbb P \left\{ \min_{i=1,\hdots,n \atop \bX_i \in \mathcal B^c_{\bx,\varepsilon}}\rho_n(\bx,\bX_i) <\gamma_n\right\}\leq 2^{d-1}n \exp \left[-n\left(\frac{V_d}{2^{2d}}\varepsilon^d f(\bx)\right)^2/2\right].$$
The upper bound goes to $0$ as $n \to \infty$.
\hfill $\blacksquare$
\paragraph{Analysis of ${\bf B}$.} Consistency of the second term in inequality (\ref{grosseinegalite}) is established in the following proposition.
\begin{pro}
\label{ulm2}
For $\mu$-almost all $\bx \in \mathbb R^d$, all $\varepsilon>0$ and all $\delta>0$ small enough, 
$$\mathbb P \left\{ \max_{i=1,\hdots,n \atop \bX_i \in \mathcal B_{\bx,\delta}}\rho_n(\bx,\bX_i) \geq \gamma_n\right\} \to 0 \quad \mbox{as } n \to \infty,$$
provided
\begin{equation}
\label{gamman}
\gamma_n=n^d\left( \frac{V_d }{2^{2d+1}}\varepsilon^d f(\bx)\right)^{2^{d-1}}.
\end{equation}
\end{pro}
{\bf Proof of Proposition \ref{ulm2}}\quad Fix $\bx$ in a set of $\mu$-measure $1$ such that $f(\bx)>0$ and denote by $N_{\bx, \delta}$ the number of hyperplanes that cut the ball $\mathcal B_{\bx, \delta}$.
Clearly,
$$
\mathbb P \left\{ \max_{i=1,\hdots,n \atop \bX_i \in \mathcal B_{\bx,\delta}}\rho_n(\bx,\bX_i) \geq  \gamma_n\right\} \leq \mathbb P \left \{ N_{\bx,\delta} \geq \gamma_n\right\}.$$
Observe that, with probability 1,
$$N_{\bx, \delta} =\sum_{1\leq i_1<\hdots <i_d\leq n} \mathbf 1_{\{\mathcal H(\bX_{i_1}, \hdots,  \bX_{i_d}) \cap \mathcal B_{\bx, \delta} \neq \emptyset\}},$$
whence, since $\bX_1, \hdots, \bX_n$ are identically distributed,
\begin{align*}
 \mathbb E[N_{\bx, \delta}] &={n \choose d} \mathbb P\left\{\mathcal H(\bX_{1}, \hdots,  \bX_{d}) \cap \mathcal B_{\bx, \delta} \neq \emptyset\right\}\\
& \leq \frac{n^d}{d!} \mathbb P\left\{\mathcal H(\bX_{1}, \hdots,  \bX_{d}) \cap \mathcal B_{\bx, \delta} \neq \emptyset\right\}.
\end{align*}
Consequently, given the choice (\ref{gamman}) for $\gamma_n$ and the result of technical Lemma \ref{vendredi4}, it follows that
$$ \mathbb E [N_{\bx,\delta}]< \gamma_n/2$$
for all $\delta$ small enough, independently of $n$. Thus, using the bounded difference inequality (McDiarmid \cite{McDiarmid}), we obtain, still with the choice
$$\gamma_n=n^d\left( \frac{V_d }{2^{2d+1}}\varepsilon^d f(\bx)\right)^{2^{d-1}},$$
\begin{align*}
\mathbb P \left \{ N_{\bx,\delta} \geq \gamma_n\right\}& \leq \mathbb P \left\{N_{\bx,\delta}-\mathbb E[N_{\bx,\delta}]\geq {\gamma_n}/{2}\right\}\\
&\leq  \exp\left(-2\frac{(\gamma_n/2)^2}{n^{2d-1}}\right)\\
&=\exp\left[-\left( \frac{V_d }{2^{2d+1}}\varepsilon^d f(\bx)\right)^{2^{d}}n/2\right].
\end{align*}
This upper bound goes to zero as $n$ tends to infinity, and this concludes the proof of the proposition.\hfill $\blacksquare$
\paragraph{Analysis of ${\bf C}$.} To achieve the proof of identity (\ref{convergenceenproba}), it remains to show that the third and last term of (\ref{grosseinegalite}) converges to $0$. This is done in the following proposition.
\begin{pro}
\label{ulm3}
Assume that $k_n/n\to 0$ as $n\to \infty$. Then, for $\mu$-almost all $\bx \in \mathbb R^d$ and all $\delta>0$,
$$\mathbb P \left \{\emph{Card}\left\{i=1,\hdots,n:\|\bX_i-\bx\| \leq \delta\right\}<k_n\right\} \to 0 \quad \mbox{as } n\to \infty. $$
\end{pro}
{\bf Proof of Proposition \ref{ulm3}}\quad Recall that the collection of all $\bx$ with $\mu(\mathcal B_{\bx,\tau})>0$ for all $\tau>0$ is called the support of $\mu$, and note that it may alternatively be defined as the smallest closed subset of $\mathbb R^d$ of $\mu$-measure 1 (Parthasarathy \cite[Chapter 2]{Parthasarathy}). Thus, fix $\bx$ in the support of $\mu$ and set
$$p_{\bx,\delta}=\mathbb P\{\bX \in \mathcal B_{\bx,\delta}\},$$
so that $p_{\bx,\delta}>0$. Then the following chain of inequalities is valid:
\begin{align*}
&\mathbb P \left\{\mbox{Card}\left\{i=1,\hdots,n:\|\bX_i-\bx\| \leq \delta\right\}<k_n\right\} \\
&\quad = \mathbb P\left\{\mbox{Binomial }(n,p_{\bx,\delta})<k_n\right\}\\
&\quad \leq \mathbb P\left\{\mbox{Binomial }(n,p_{\bx,\delta})\leq {np_{\bx,\delta}}/{2}\right\}\\
&\qquad (\mbox{for all $n$ large enough, since $k_n/n$ tends to $0$})\\
&\quad \leq \exp (-np^2_{\bx,\delta}/2),
\end{align*}
where the last inequality follows from Hoeffding's inequality (Hoeffding \cite{Hoeffding}). This terminates the proof of Proposition \ref{ulm3}.
\hfill $\blacksquare$
\bigskip

We have proved so far that, for $\mu$-almost all $\bx$, as $k_n/n\to0$, the quantity $\max_{i=1,\hdots,k_n}\| \bX_{(i)}(\bx)-\bx\|$ converges to $0$ in probability. By the elementary inequality 
$$\mathbb E \left [\frac{1}{k_n}\sum_{i=1}^{k_n}\mathbf 1_{\{\|\bX_{(i)}(\bx)-\bx\|>\varepsilon\}}\right ]\leq \mathbb P \left\{ \max_{i=1, \hdots, k_n} \|\bX_{(i)}(\bx)-\bx\|>\varepsilon \right\},$$
it immediately follows that, for such an $\bx$,
\begin{equation}
\label{printemps}
\mathbb E \left[\frac{1}{k_n}\sum_{i=1}^{k_n}\mathbf 1_{\{\|\bX_{(i)}(\bx)-\bx\|>\varepsilon\}}\right]\to 0
\end{equation}
provided $k_n/n \to 0$.
We are now ready to complete the proof of Theorem \ref{theoreme1}.
\medskip

Fix $\bx$ in a set of $\mu$-measure 1 such that consistency (\ref{printemps}) holds and $r$ is continuous at $\bx$ (this is possible by the assumption on $r$). Because $|a+b|^p \leq 2^{p-1}(|a|^p +|b|^p)$ for $p\geq 1$, we see that
\begin{align*}
\mathbb E \left|r_n(\bx)-r(
\bx)\right|^p &\leq 2^{p-1} \mathbb E \left| \frac{1}{k_n} \sum_{i=1}^{k_n} \left[Y_{(i)}(\bx) -r\left(
\bX_{(i)}(\bx)\right)\right]\right|^p\\
& \quad +2^{p-1} \mathbb E \left|\frac{1}{k_n} \sum_{i=1}^{k_n} \left[r\left(
\bX_{(i)}(\bx)\right)-r(\bx)\right]\right|^p.
\end{align*}
Thus,  by Jensen's inequality,
\begin{align*}
\mathbb E \left|r_n(\bx)-r(
\bx)\right|^p &\leq 2^{p-1} \mathbb E \left| \frac{1}{k_n}\sum_{i=1}^{k_n} \left[Y_{(i)}(\bx) -r\left(
\bX_{(i)}(\bx)\right)\right]\right|^p\\
& \quad +2^{p-1} \mathbb E \left[\frac{1}{k_n}\sum_{i=1}^{k_n} \left|r\left(
\bX_{(i)}(\bx)\right)-r(
\bx)\right|^p\right]\\
&  :=2^{p-1}{\mathbf I_n}+2^{p-1}{\mathbf J_n}. 
\end{align*}
Firstly, for arbitrary $\varepsilon >0$, we have
\begin{align*}
{\mathbf J_n}&= \mathbb E \left[Ê\frac{1}{k_n}\sum_{i=1}^{k_n} \left |r\left(
\bX_{(i)}(\bx)\right) -r(\bx)\right|^p\mathbf 1_{\{\|\bX_{(i)}(\bx)-\bx\|>\varepsilon\}} \right] \\
& \quad +  \mathbb E \left[Ê\frac{1}{k_n}\sum_{i=1}^{k_n} \left |r\left(
\bX_{(i)}(\bx)\right) -r(
\bx)\right|^p\mathbf 1_{\{\|\bX_{(i)}(\bx)-\bx\|\leq\varepsilon\}}\right],
\end{align*}
whence
\begin{align*}
{\mathbf J_n}&\leq 2^p\zeta^p\,\mathbb E \left[Ê\frac{1}{k_n}\sum_{i=1}^{k_n}\mathbf 1_{\{\|\bX_{(i)}(\bx)-\bx\|>\varepsilon\}} \right] \\
& \quad +\left[\sup_{\mathbf y \in \mathbb R^d : \|\mathbf y-\bx\|\leq \varepsilon} \left| r(
\mathbf y)-r(\bx)\right|\right]^p\\
& \qquad (\mbox{since $|Y|\leq \zeta$}).
\end{align*}
The first term on the right-hand side of the latter inequality tends to 0 by (\ref{printemps}) as $k_n/n \to 0$, whereas the rightmost one can be made arbitrarily small as $\varepsilon \to 0$ since $r$ is continuous at $\bx$. This proves that ${\mathbf J_n} \to 0$ as $n \to \infty$.
\medskip

Next, by successive applications of inequalities of Marcinkiewicz and Zygmund \cite{Marcin} (see also Petrov \cite[pages 59-60]{Petrov}), we have for some positive constant $C_p$ depending only on $p$,
\begin{align*}
& {\mathbf I_n}\leq C_p\,\mathbb E \left[ \frac{1}{k_n^2} \sum_{i=1}^{k_n} \left|Y_{(i)}(\bx) -r\left(
\bX_{(i)}(\bx)\right)\right|^2\right]^{p/2} \\
 & \quad \leq \frac{(2\zeta)^{p}C_p}{k_n^{p/2}}\\
 & \qquad (\mbox{since $|Y|\leq \zeta$}).
\end{align*}
Consequently, ${\mathbf I_n} \to 0$ as $k_n\to \infty$, and this concludes the proof of the theorem.
\section{Some technical lemmas}
The notation of this section is identical to that of Section 3. In particular, it is assumed throughout that $\bX$ has a probability density $f$ with respect to the Lebesgue measure $\lambda$ on $\mathbb R^d$. This requirement implies that any collection $\bX_{i_1},\hdots,\bX_{i_d}$ ($1\leq i_1<i_2<\hdots<i_d\leq n)$ of $d$ points among $\bX_1, \hdots, \bX_n$ define with probability 1 a unique hyperplane $\mathcal H(\bX_{i_1},\hdots,\bX_{i_d})$ in $\mathbb R^d$. Recall finally that, for $\bx \in \mathbb R^d$ and $\varepsilon >0$, we set
$$p_{\bx, \varepsilon}=\min_{j=1, \hdots, 2^{d-1}}\inf_{{\Theta}} \mu \left\{ \mathcal R^j_{\bx,\varepsilon}({\Theta})\right\},$$
where the infimum is taken over all possible hyperspherical angles ${\Theta}$, and the regions $\mathcal R^j_{\bx,\varepsilon}({\Theta})$, $j=1, \hdots, 2^{d-1}$, define a partition of the ball $\mathcal B_{\bx,\varepsilon}(\Theta)$. Recall also that the numbers of sample points falling in each of these regions are denoted by $N^1_{\bx,\varepsilon}({\Theta}), \hdots, N^{2^{d-1}}_{\bx,\varepsilon}({\Theta})$. For a better understanding of the next lemmas, the reader should refer to Figure \ref{angle} and Figure \ref{anglebis}.
\begin{lem}
\label{samedi26}
For $\mu$-almost all $\bx \in \mathbb R^d$ and all $\varepsilon>0$ small enough,
 $$p_{\bx, \varepsilon}\geq \frac{V_d}{2^{2d}}{\varepsilon}^df(\bx)>0.$$
\end{lem}
{\bf Proof of Lemma \ref{samedi26}}\quad 
We let $\bx$ be a Lebesgue point of $f$, that is, an $\bx$ such that for any collection $\mathcal A$ of subsets of $\mathcal B_{\mathbf 0,1}$ with the property that for all $A \in \mathcal A$, $\lambda(A)\geq c \lambda(\mathcal B_{\mathbf 0,1})$ for some fixed $c>0$,
\begin{equation}
\label{lebesgue}
\lim_{\varepsilon \to 0} \sup_{A \in \mathcal A} \left| \frac{\displaystyle \int_{\bx+\varepsilon A} f(\mathbf y)\mbox{d}\mathbf y}{\lambda\{\bx+\varepsilon A\}}-f(\bx)\right|=0,
\end{equation}
where $\bx+\varepsilon A=\{\mathbf y \in \mathbb R^d:(\mathbf y-\bx)/\varepsilon \in A\}$. As $f$ is a density, we know that $\mu$-almost all $\bx$ satisfy this property (see, for instance, Wheeden and Zygmund \cite{Wheeden}). Moreover, since $f$ is $\mu$-almost surely positive, we may also assume that $f(\bx)>0$.
\medskip

Thus, keep such an $\bx$ fixed. Fix also $j \in \{1, \hdots, 2^{d-1}\}$, and set
$$p^j_{\bx, \varepsilon}=\inf_{{\Theta}} \mu \left\{ \mathcal R^j_{\bx,\varepsilon}({\Theta})\right\}.$$
Taking for $\mathcal A$ the collection of regions $\mathcal R^j_{\mathbf 0,1}({\Theta})$ when the hyperspherical angle ${\Theta}$ varies, that is,
$$\mathcal A=\left\{ \mathcal R^j_{\mathbf 0,1}({\Theta}) : {\Theta} \in  [0,\pi]^{d-2} \times [0,2\pi)\right\},$$
and observing that
$$\lambda \left \{\mathcal R^j_{\bx,\varepsilon}({\Theta})\right\}=\frac{V_d}{2^{d-1}}\left (\frac{\varepsilon}{2}\right)^d,$$
we may write, for each $j=1, \hdots, 2^{d-1}$,
\begin{align*}
\left |\frac{2^{d-1}p^j_{\bx, \varepsilon}}{V_d(\varepsilon/2)^d}-f(\bx)\right| & = \left | \inf_{{\Theta}} \frac{\mu \left\{ \mathcal R^j_{\bx,\varepsilon }({\Theta})\right\}}{\lambda \left \{ \mathcal R^j_{\bx,\varepsilon }({\Theta})\right\}} -f(\bx)\right|\\
& =\left|\inf_{A \in \mathcal A} \frac{\displaystyle \int_{\bx+\varepsilon A} f(\mathbf y)\mbox{d}\mathbf y}{\lambda\{\bx+\varepsilon A\}}-f(\bx)\right|\\
& \leq \sup_{A \in \mathcal A} \left| \frac{\displaystyle \int_{\bx+\varepsilon A} f(\mathbf y)\mbox{d}\mathbf y}{\lambda\{\bx+\varepsilon A\}}-f(\bx)\right|.
\end{align*}
The conclusion follows from identity (\ref{lebesgue}).
\hfill $\blacksquare$
\begin{lem}
\label{lundi28}
Fix $\bx \in  \mathbb R^d$, $\varepsilon >0$ and ${\Theta} \in [0,\pi]^{d-2} \times [0,2\pi)$. Let $N_{\bx,\varepsilon}({\Theta})$ be the number of hyperplanes passing through $d$ out of the observations $\bX_1,\hdots,\bX_n$ and cutting the segment $\mathcal S_{\bx,\varepsilon}({\Theta})$. Then, with probability 1,
$$N_{\bx,\varepsilon}({\Theta}) \geq \frac{N^1_{\bx,\varepsilon}({\Theta}) \hdots N^{2^{d-1}}_{\bx,\varepsilon}({\Theta})}{n^{2^{d-1}-d}}.$$
\end{lem}
{\bf Proof of Lemma \ref{lundi28}}\quad If one of the $N_{\bx,\varepsilon}^j(\Theta)$ $(j=1, \hdots, 2^{d-1})$ is zero, then the result is trivial. Thus, in the rest of the proof, we suppose that each $N_{\bx,\varepsilon}^j(\Theta)$ is positive and note that this implies $n \geq 2^{d-1}$.
\medskip

Pick sequentially $2^{d-1}$ observations, say $\bX_{i_1}, \hdots, \bX_{i_{2^{d-1}}}$, in the $2^{d-1}$ regions $\mathcal R^1_{\bx, \varepsilon}({\Theta}), \hdots, \mathcal R^{2^{d-1}}_{\bx, \varepsilon}({\Theta})$. By construction, the polytope defined by these $2^{d-1}$ points cuts the axe $\mathcal L_{\bx,\varepsilon}(\Theta)$. Consequently, with probability 1, any hyperplane drawn according to $d$ out of these $2^{d-1}$ points cuts the segment $\mathcal S_{\bx,\varepsilon}({\Theta})$. The result follows by observing that there are exactly $N^1_{\bx, \varepsilon}({\Theta})\hdots N^{2^{d-1}}_{\bx, \varepsilon}({\Theta})$ such polytopes.
\hfill $\blacksquare$
\begin{lem}
\label{vendredi4}
For $1\leq i_1 <\hdots <i_d\leq n$, let $\mathcal H(\bX_{i_1},\hdots,\bX_{i_d})$ be the hyperplane passing through $d$ out of the observations $\bX_1,\hdots,\bX_n$. Then, for all $\bx \in \mathbb R^d$,
$$\mathbb P \left \{ \mathcal H(\bX_{i_1}, \hdots, \bX_{i_d}) \cap \mathcal B_{\bx, \delta} \neq \emptyset \right\}\to 0 \quad \mbox{as } \delta \downarrow 0.$$
\end{lem}
{\bf Proof of Lemma \ref{vendredi4}}\quad Given two hyperplanes $\mathcal H$ and $\mathcal H'$ in $\mathbb R^d$, we denote by $\Phi(\mathcal H, \mathcal H')$ the (dihedral) angle between $\mathcal H$ and $\mathcal H'$. Recall that $\Phi(\mathcal H, \mathcal H') \in [0,\pi]$ and that it is defined as the angle between the corresponding normal vectors.
\medskip

\begin{figure}[h]
\psfrag{x}{$\bX_{i_1}$}
\psfrag{y}{$\bX_{i_2}$}
\psfrag{d}{$\delta$}
\psfrag{c}{$\bx$}
\psfrag{t}{$\Phi_{\delta}$}
\centering
\includegraphics*[height=8cm]{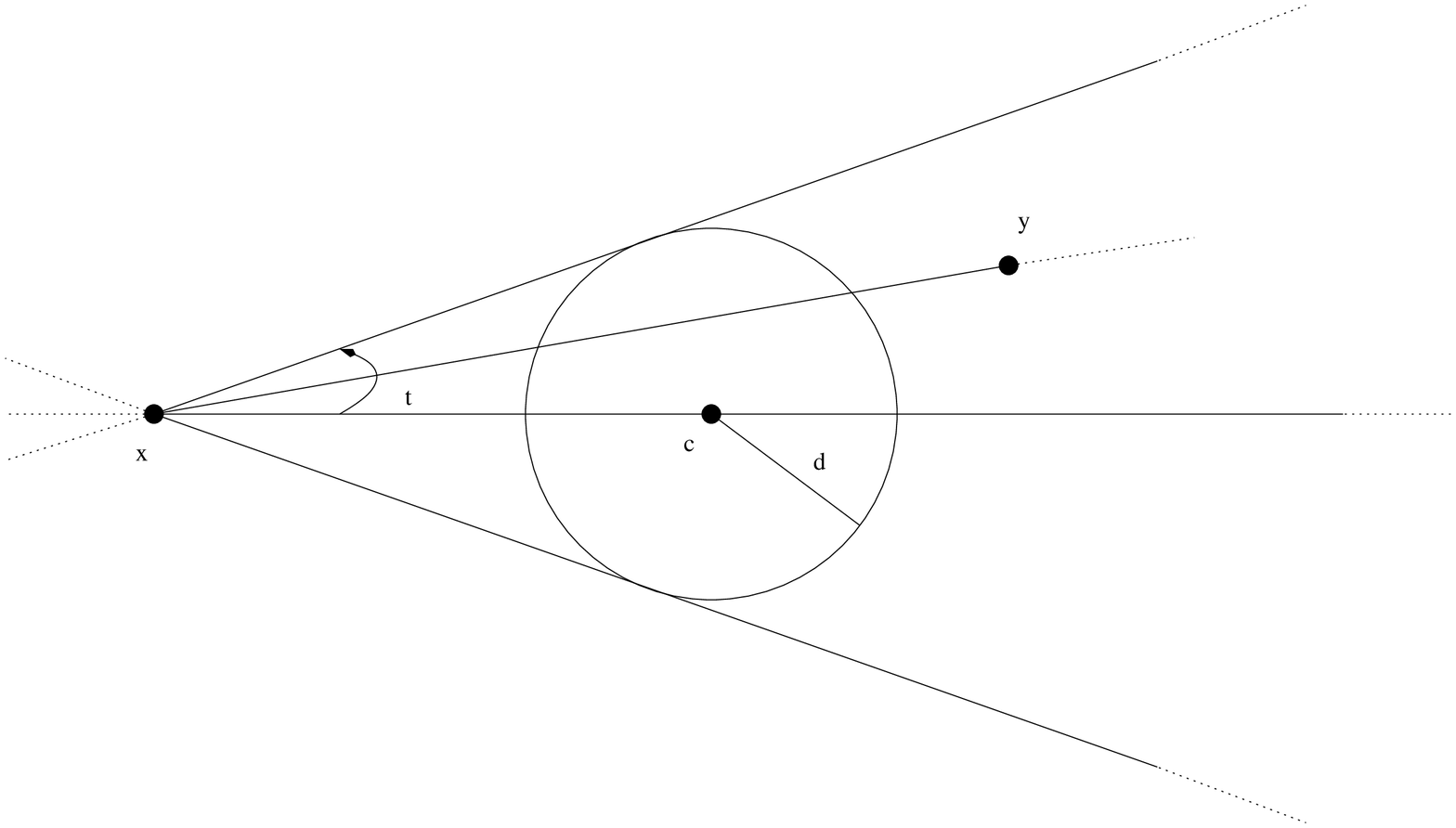}
\caption{The hyperplanes $\mathcal H(\bx, \bX_{i_1}, \hdots, \bX_{i_{d-1}})$ and $\mathcal H(\bX_{i_1}, \hdots, \bX_{i_{d}})$, and the angle $\Phi_{\delta}$. Illustration in dimension 2.}
\label{tangent}
\end{figure}

Fix $1\leq i_1 <\hdots <i_d\leq n$. Let $\mathcal E_{\delta}$ be the event
$$\mathcal E_{\delta}=\left \{ \|\bX_{i_j}-\bx\| >\delta : j=1, \hdots, d-1\right\},$$
and let $\mathcal H( \bx, \bX_{i_1}, \hdots, \bX_{i_{d-1}})$ be the hyperplane passing through $\bx$ and the $d-1$ points $\bX_{i_1}, \hdots, \bX_{i_{d-1}}$. Clearly, on $\mathcal E_{\delta}$, the event $\{\mathcal H(\bX_{i_1}, \hdots, \bX_{i_{d}}) \cap \mathcal B_{\bx, \delta} \neq \emptyset\}$ is the same as
$$\left \{\Phi\left (\mathcal H(\bx, \bX_{i_1}, \hdots, \bX_{i_{d-1}}) , \mathcal H(\bX_{i_1}, \hdots, \bX_{i_{d}})\right) \leq \Phi_{\delta}\right\},$$
where $\Phi_{\delta}$ is the angle formed by $\mathcal H(\bx, \bX_{i_1}, \hdots, \bX_{i_{d-1}})$ and the hyperplane going trough $\bX_{i_1}, \hdots, \bX_{i_{d-1}}$ and tangent to $\mathcal B_{\bx, \delta}$ (see Figure \ref{tangent} for an example in dimension 2). 
\medskip

Thus, with this notation, we may write
\begin{align*}
& \mathbb P \left \{ \mathcal H(\bX_{i_1}, \hdots, \bX_{i_d}) \cap \mathcal B_{\bx, \delta} \neq \emptyset \right\} \\
& \quad \leq \mathbb P\{\mathcal E_n^c\} + \mathbb P \left \{ \mathcal H(\bX_{i_1}, \hdots, \bX_{i_d}) \cap \mathcal B_{\bx, \delta} \neq \emptyset, \mathcal E_n \right\}\\
& \quad \leq  \mathbb P\{\mathcal E_n^c\} + \mathbb P \left \{\Phi\left (\mathcal H(\bx, \bX_{i_1}, \hdots, \bX_{i_{d-1}}) , \mathcal H(\bX_{i_1}, \hdots, \bX_{i_{d}})\right) \leq \Phi_{\delta},\mathcal E_{n} \right\}.
\end{align*}
Since $\bX$ has a density, the first of the two terms above tends to zero as $\delta \downarrow 0$. To analyze the second term, first note that, conditionally on $\bX_{i_1}, \hdots, \bX_{i_{d-1}}$, the angle $\Phi (\mathcal H(\bx, \bX_{i_1}, \hdots, \bX_{i_{d-1}}) , \mathcal H(\bX_{i_1}, \hdots, \bX_{i_{d}}))$ is absolutely continuous with respect to the Lebesgue measure on $\mathbb R$. This follows from the following two observations: $i)$ the random variable $\bX_{i_d}$ has a density with respect to the Lebesgue measure on $\mathbb R^d$, and $ii)$ conditionally on $\bX_{i_1}, \hdots, \bX_{i_{d-1}}$, $\Phi (\mathcal H(\bx, \bX_{i_1}, \hdots, \bX_{i_{d-1}}) , \mathcal H(\bX_{i_1}, \hdots, \bX_{i_{d}}))$ is obtained from $\bX_{i_d}$ via translations, orthogonal transformations, and the arctan function. 
\medskip

Thus, writing
\begin{align*}
& \mathbb P \left \{\Phi\left (\mathcal H(\bx, \bX_{i_1}, \hdots, \bX_{i_{d-1}}) , \mathcal H(\bX_{i_1}, \hdots, \bX_{i_{d}})\right) \leq \Phi_{\delta} ,\mathcal E_n\right\} \\
& \!= \mathbb E \!\left [\mathbf 1_{\mathcal E_n}\mathbb P \left \{\Phi\left (\mathcal H(\bx, \bX_{i_1}, \hdots, \bX_{i_{d-1}}) , \mathcal H(\bX_{i_1}, \hdots, \bX_{i_{d}})\right) \leq \Phi_{\delta} | \bX_{i_1}, \hdots, \bX_{i_{d-1}}\right\}\right]
\end{align*}
and noting that, on the event $\mathcal E_n$, for fixed $\bX_{i_1}, \hdots, \bX_{i_{d-1}}$, $\Phi_{\delta} \downarrow 0$ as $\delta \downarrow 0$, we conclude by the Lebesgue dominated convergence theorem that
$$\mathbb P \left \{\Phi\left (\mathcal H(\bx, \bX_{i_1}, \hdots, \bX_{i_{d-1}}) , \mathcal H(\bX_{i_1}, \hdots, \bX_{i_{d}})\right) \leq \Phi_{\delta} \right\} \to 0\quad \mbox{as }  \delta \downarrow 0.$$
\hfill $\blacksquare$
\paragraph{Acknowledgments.} We thank two anonymous referees for valuable comments and insightful suggestions. 
\bibliographystyle{plain}
\bibliography{biblio-ppv}
\end{document}